\documentstyle[12pt,dina4]{article}
\newcommand{\ve}{{\varepsilon}}

\newtheorem{rem}{Remark}[section]

\newtheorem{LEMMA}{Lemma}

\newtheorem{COROLLARY}{Corollary}
\newcommand{\coro}[1]{
\begin{COROLLARY}
#1
\end{COROLLARY}
}
\newtheorem{THEOREM}{Theorem}[section]
\newcommand{\thm}[1]{
\begin{THEOREM}
#1
\end{THEOREM}
}
\newtheorem{DEFINITION}{Definition}[section]

\newcommand{\pf}{
\begin{flushleft}
{\bf Proof}
\end{flushleft}}
\newcommand{\nm}[1]{\parallel #1 \parallel}

\newcommand{\nMw}[1]{\nm{#1}_{M^p_{w(\alpha )}}}

\newcommand{\nmwg}[1]{\nm{#1}_{L^p_{w(\gamma )}}}
\newcommand{\nmgk}[1]{\nm{#1}_{L^1_{w(\gamma )}}}
\newcommand{\nMwg}[1]{\nm{#1}_{M^p_{w(\gamma )}}}

\begin{document}
 \begin{center}
 {\LARGE{}On necessary multiplier conditions for
 Laguerre expansions II}\\[.4cm]
 {\sc{}George  Gasper\footnote{Department of Mathematics,
Northwestern
 University, Evanston, IL 60208, USA. The work of this author was
 supported in part by the National Science   Foundation under grant
 DMS-9103177.}
 and Walter  Trebels\footnote{Fachbereich Mathematik, TH Darmstadt,
D--6100
 Darmstadt, Germany.} }\\[.4cm]
 {\it Dedicated to Dick Askey and Frank Olver}\\[.4cm]
 {(March 5, 1992 version)}\\[.4cm]
 \end{center}
 
 \bigskip
 {\bf Abstract.}
The necessary multiplier conditions for Laguerre expansions
derived in Gasper and
Trebels \cite{laguerre} are supplemented and modified. This allows
us
to place Markett's Cohen type inequality \cite{cohen} (up to the
$\log $--case) in the general framework of necessary conditions.
 
 
\bigskip
{\bf AMS(MOS) subject classifications.} 33C65, 42A45, 42C10
 
\bigskip
\section{ Introduction}
 
The purpose of this sequel to \cite{laguerre} is
to obtain a better insight into the structure of
Laguerre multipliers on $L^p$ spaces from the point of view of
necessary
conditions. We recall that in \cite{laguerre} there occurs the 
annoying
phenomenon that, e.g., the optimal necessary conditions in the case
$p=1$ do
not give the ``right'' unboundedness behavior of the Ces\`aro
means. By
slightly modifying these conditions we can not only remedy this
defect but can
also derive Markett's Cohen type inequality \cite{cohen} (up to the
$\log$--case) as an immediate consequence.\\
For the convenience of the reader we briefly repeat the notation;
we consider the  Lebesgue spaces
$$L^p_{w(\gamma )} = \{ f: \; \nmwg{f} = ( \int _0^\infty
|f(x) e^{-x/2}|^p x^{\gamma }\, dx )^{1/p} < \infty \} \; ,\quad 1
\le p <
\infty  ,$$
denote the classical Laguerre polynomials
by $L_n^\alpha (x), \; \alpha >-1,\; n\in {\bf N}_0$
  (see Szeg\"o \cite[p. 100]{szego}), and set
$$ R_n^\alpha (x)=L_n^\alpha (x)/L_n^\alpha (0),\; \; \; 
L_n^\alpha (0)=
A_n^\alpha =
\left( \matrix{ n+\alpha \cr n} \right) =\frac{\Gamma (n+ \alpha
+1)}{\Gamma
(n+1) \Gamma (\alpha +1)  }.$$
Associate to $f$ its formal Laguerre series
$$ f(x) \sim (\Gamma (\alpha +1))^{-1}
\sum _{k=0}^\infty {\hat f}_\alpha (k) L_k^\alpha (x), $$
where the Fourier Laguerre coefficients of $f$ are defined by
\begin{equation}\label{coeff}
{\hat f}_\alpha (n) = \int _0^\infty f(x)R_n^\alpha (x)x^\alpha
e^{-x}\, dx
\end{equation}
(if the integrals exist).
A sequence $m=\{ m_k\} $ is called a (bounded) multiplier on
$L^p_{w(\gamma )}$, notation $ m \in M^p_{w(\gamma )}$, if
$$ \nmwg{
\sum _{k=0}^\infty m_k{\hat f}_\alpha (k) L_k^\alpha }\le C
\nmwg{f}$$
for all polynomials $f$; the smallest constant $C$ for which this
holds is
called
the multiplier norm $ \nMwg{ m }. $ The necessary conditions will
be
given in certain ``smooth\-ness'' properties of the multiplier
sequence in
question. To this end we introduce a
 fractional difference operator of order $\delta$ by
$$ \Delta ^\delta m_k = \sum _{j=0}^\infty A_j^{-\delta
-1}m_{k+j}$$
(whenever the sum converges), the first order difference operator
$\Delta_2$ with increment $2$ by
$$\Delta_2 m_k=m_k-m_{k+2},$$
and the notation
$$ \Delta_2 \Delta ^\delta m_k = \Delta ^{\delta +1} m_k +
\Delta ^{\delta +1}m_{k+1}. $$
Generic positive constants that are independent of the functions
(and sequences) will be denoted by $C$.
Within the setting of the $L^p_{w(\gamma )}$-spaces our main
results now read
(with $1/p+1/q=1$):
\thm{
Let $ \alpha ,\, a>-1$ and $\alpha +a>-1$.
If  \mbox{$f\in L^p_{w(\gamma )},\; 1\le p <2,$}
 then
\begin{equation}\label{fdiff}
 \left( \sum _{k=0}^\infty |(k+1)^{(\gamma +1)/p - 1/2} \Delta _2
\Delta^a
{\hat f}_\alpha (k)|^q \right) ^{1/q} \le C \nmwg{f} ,
\end{equation}
provided
 $$\frac{\gamma +1}{p} \le \frac{\alpha +a}{p}+1 \quad \quad \quad
\quad \quad
{\it if} \; \,
\alpha +a\le 1/2, $$
$$\frac{\gamma +1}{p} \le \frac{\alpha +
a}{2}+1+\frac{1}{2}\left( \frac{1}{p} -\frac{1}{2}\right)  \quad
{\it if} \; \,
\alpha +a > 1/2. $$
}
 
\bigskip \noindent
As in \cite{laguerre} (see there the proof of Lemma 2.3)
we immediately obtain
\thm{
Let $m=\{ m_k\} \in M^p_{w(\gamma )},\, 1\le p <2,$ and let $\alpha
$ and $a$
be as in Theorem 1.1. Then
\begin{equation}\label{mdiff}
\sup _n \left( \sum _{k=n}^{2n}
|(k+1)^{(2\gamma +1)/p -(2\alpha + 1)/2} \Delta _2 \Delta^{a}
m_k|^q\frac{1}{k+1}\right) ^{1/q} \le C \nMwg{ m },
\end{equation}
provided that in the case $\alpha +a\le 1/2$ the condition
$$\frac{\alpha +a}{p}+1 \ge  \frac{\gamma +1}{p} >
\cases{
(\alpha  +1)/2 +1/3p & if $\; 1\le p<4/3$ \cr
(\alpha +1)/2 +1/4  &  if $\; 4/3 \le p <2$ }
$$
holds, and in the case $\alpha +a >1/2$ the condition
$$\frac{\alpha +a}{2}+1+\frac{1}{2}\left( \frac{1}{p}
-\frac{1}{2}\right)
\ge \frac{\gamma +1}{p} >
\cases{
(\alpha  +1)/2 +1/3p & if $\; 1\le p<4/3$ \cr
(\alpha +1)/2 +1/4  &  if $\; 4/3 \le p <2$. }
 $$
}
 
\bigskip \noindent
In view of the results in \cite{cohen}, \cite{laguerre} and
for an easy comparison we want to emphasize the cases
$\gamma = \alpha$ and $\gamma = \alpha p/2$. Therefore, we state
\coro{
\begin{itemize}
\item[a)] Let $m\in M^p_{w(\alpha )},\; 1\le p < 2,$ and let $
\alpha >-1$
be such that $\max \{ 1/(3p),\, 1/4\} < (\alpha +1)(1/p-1/2)$.
Then,
with $\lambda : = (2\alpha +1)(1/p-1/2) $,
$$ \sup _n \left( \sum _{k=n}^{2n}
|(k+1)^\lambda \Delta _2 \Delta^{\lambda -1} m_k|^q
\frac{1}{k+1} \right) ^{1/q} \le C \nMw{ m }.$$
\item[b)] Let $m\in M^p_{w(\alpha p/2)},\; 1\le p < 4/3$, and $
(\alpha
-1)(1/p-1/2) \ge -1/2$ .
Then
$$ \sup _n \left( \sum _{k=n}^{2n}
|(k+1)^{1/p-1/2}\Delta _2 m_k|^q
\frac{1}{k+1} \right) ^{1/q} \le C \nm{ m }_{M^p_{w(\alpha
p/2)}}.$$
\end{itemize}
}
 
\bigskip \noindent
{\bf Remarks.}
1) For polynomial $f(x)=\sum_{k=0}^nc_k L_k^\alpha (x)$
Theorem 1.1 yields, by taking only the $(k=n)$--term on the left
hand side of (\ref{fdiff}),
$$ |c_n|(n+1)^{(\gamma +1)/p-1/2} \le C \nmwg{f}, \quad 1\le p<2 $$
(under the restrictions on $\gamma $ of Theorem 1.1).
In particular, if we choose $\gamma =\alpha$,
this comprises  formula (1.13) in Markett \cite{cohen} for his
basic
case $\beta =\alpha$.
For $\gamma =\alpha p/2$, it even extends
formula (1.14) in \cite{cohen} to negative $\alpha $'s as described
in
Corollary 1.3,~b).  The case
$2<p<\infty $ can be done by an application of a Nikolskii
inequality,
see \cite{cohen}.
 
\bigskip \noindent
2) Analogously, Cohen type inequalities follow from Theorem 1.2; in
particular,
Corollary 1.3 yields
\coro{
Let $m=\{ m_k \} _{k=0}^n $ be a finite sequence, $1\le p <2$, and
$\alpha >-1$.
\begin{itemize}
\item[a)] If $m\in M^p_{w(\alpha )}$ then
$$ (n+1)^{(2\alpha +2)(1/p-1/2)-1/2}|m_n|\le C \nMw{m}, \quad 1\le
p<
\frac{4\alpha +4}{2\alpha +3},$$
provided $\max\{ 1/3p, \, 1/4 \} < (\alpha +1)(1/p-1/2)$.
\item[b)] If $m\in M^p_{w(\alpha p/2)}$ and $(\alpha -1)(1/p-1/2)
\ge -1/2$,
 then
$$ (n+1)^{2/p -3/2}|m_n|\le C \nm{m} _{M^p_{w(\alpha p/2)}}, \quad
1\le p< 4/3.
$$
\end{itemize}
}
\bigskip
With the exception of the crucial $\log$--case, i.e. $p_0=(4\alpha
+4)/(2\alpha +3)$ or $p_0=4/3$, resp., Corollary 1.4 contains
Markett's Theorem
1
in \cite{cohen} and extends it to negative $\alpha$'s.
In particular we obtain for the Ces\`aro means
of order $\delta \ge 0$, represented by its multiplier
sequence $m^\delta_{k,n}=A_{n-k}^\delta /A_n^\delta$, the ``right''
unboundedness behavior (see \cite{indag} )
$$\nMw{\{ m^\delta_{k,n} \} } \ge C (n+1)^{(2\alpha
+2)(1/p-1/2)-1/2 - \delta
}, \quad 1\le p <\frac{4\alpha +4}{2\alpha +3+2\delta}.$$
 
\medskip \noindent
3) There arises the question, in how far the type of necessary
conditions in
\cite{laguerre} are comparable with the present ones. Let $\lambda
>1$.
 Since $\Delta_2m_k= \Delta m_k + \Delta m_{k+1}$ we  obviously
have
\begin{equation}\label{embedding}
\sup _n \left( \sum _{k=n}^{2n}
|(k+1)^{(2\gamma +1)/p -(2\alpha + 1)/2} \Delta _2 \Delta ^{\lambda
-1}
m_k|^q\frac{1}{k+1}\right) ^{1/q} \quad \quad \quad
\end{equation}
$$ \quad \quad \quad \le C \sup _n \left( \sum _{k=n}^{2n}
|(k+1)^{(2\gamma +1)/p -(2\alpha + 1)/2} \Delta ^\lambda
m_k|^q\frac{1}{k+1}
\right) ^{1/q}.$$
In general, a converse cannot hold as can be seen by the following
example:
choose $\gamma =\alpha ,\; \lambda =(2\alpha +1)(1/p-1/2)$ and
 $m_k=(-1)^k k^{-\ve}, \; 0 < \ve < 1$. Then
$$\sup _n \left( \sum _{k=n}^{2n}
|(k+1) \Delta m_k|^q\frac{1}{k+1}\right) ^{1/q}=\infty $$
and hence by the embedding properties of the $wbv$--spaces, see
\cite{wbv}, the
right hand side of (\ref{embedding}) cannot be finite for all
$\lambda >
1$. But since $\Delta_2\Delta^{\lambda -1}m_k=
\Delta^{\lambda -1}\Delta_2 m_k \sim (k+1)^{-\ve -\lambda}$, the
left hand side
of (\ref{embedding}) is finite for all $\lambda >1$.
 
\bigskip \noindent
Theorem 1.1 will be proved in Section 2 by interpolating between
$(L^1,l^\infty )$--  and $(L^2,l^2)$--estimates. The $a\neq 0$ case
is an easy consequence of the case $a=0$
when one uses the basic formula (see formula (3) in
\cite{laguerre} and Remark 3 preceding Section 3 there)
\begin{equation}\label{gebr}
\Delta ^a R_k^\alpha (x) = \frac{\Gamma (\alpha +1)}{\Gamma (\alpha
+a+1)} x^a R_k^{\alpha +a}(x), \quad
x>0,\; a > -1 - \min \{ \alpha, \alpha /2 -1/4\},
\end{equation}
where in the case $a > -(2\alpha +1)/4$
the series for the fractional difference converges absolutely.
In Section 3, a necessary $(L^1,l^1)$--estimate is derived and it
is
compared with a corresponding sufficient $(l^1,L^1)$--estimate.
 
\bigskip
\section{ Proof of Theorem 1.1}
 
Let us first handle the $(L^2,l^2)$--estimate. Since
$$\Delta _2 \Delta^{a}{\hat f}_\alpha (k)=
\Delta^{1+a}{\hat f}_\alpha (k)+ \Delta^{1+a}{\hat f}_\alpha
(k+1)$$
it follows from the Parseval formula preceding Corollary 2.5 in
\cite{laguerre}
that
\begin{equation}\label{parseval}
\left( \sum _{k=0}^\infty |\sqrt{A_k^{\alpha +1+a}}
\Delta _2 \Delta^{a}{\hat f}_\alpha (k)|^2\right) ^{1/2}
\le C \left( \int_0^\infty |f(t)e^{-t/2}t^{(\alpha +1+a)/2}|^2
\, dt\right) ^{1/2}.
\end{equation}
Concerning the $(L^1,l^\infty )$--estimate we first restrict
ourselves to
 the case $a =0$. Define $\mu \in {\bf R}$ by
$$2\left(\frac{1}{p}-\frac{1}{2}\right) \mu =
\frac{\gamma}{p}-\frac{\alpha +1}{2} ;$$
with the notation ${\cal L}_k^\alpha (t)= (A_k^\alpha
/\Gamma(\alpha +1))^{1/2}
R_k^\alpha (t)e^{-t/2}t^{\alpha /2}$ it follows that
$$|\Delta_2 {\hat f}_\alpha (k)|=C |\int_0^\infty f(t)\{ {\cal
L}_k^\alpha (t)
/ \sqrt{A_k^\alpha} - {\cal L}_{k+2}^\alpha
(t)/\sqrt{A_{k+2}^\alpha}\}
e^{-t/2} t^{\alpha /2}\, dt| \quad \quad \quad $$
$$\le C (k+1)^{-1-\alpha/2}\int_0^\infty |f(t)|
|t^{-\mu -1/2}{\cal L}_k^\alpha (t)| e^{-t/2}
t^{(\alpha +1)/2+\mu }\, dt  \quad \quad \quad \quad \quad \quad
\quad $$
$$ +C (k+1)^{-\alpha/2}\int_0^\infty |f(t)|
|t^{-\mu -1/2}\{ {\cal L}_k^\alpha (t)
-{\cal L}_{k+2}^\alpha (t)\} | e^{-t/2} t^{(\alpha +1)/2+\mu }\, dt
= I + II.$$
We distinguish the two cases $\alpha \le 1/2$ and $\alpha >1/2$:
 
\medskip \noindent
First consider  the case $\alpha \le 1/2$ .
By the asymptotic estimates for ${\cal L}_k^\alpha (t)-{\cal
L}_{k+2}^\alpha
(t)$ in Askey and Wainger \cite[p.699]{aw}, see formula (2.12) in
\cite{cohen}, it follows for $\gamma \le \alpha +p-1$
that
$$ \nm{t^{-\mu -1/2}\{ {\cal L}_k^\alpha (t)-{\cal L}_{k+2}^\alpha
(t)\}
}_\infty \le C (k+1)^{-1-\mu}$$
so that
\begin{equation}\label{II}
 II \le C (k+1)^{-1-\mu-\alpha /2}\int_0^\infty |f(t)|e^{-t/2}
t^{(\alpha +1)/2+\mu}\, dt, \quad \gamma \le \alpha +p-1  .
\end{equation}
By Lemma 1, 4th case, in \cite{analy}
$$ \nm{t^{-\mu -1/2} {\cal L}_k^\alpha (t)}_\infty
\le C (k+1)^{-\mu -5/6}$$
so that trivially
$$ I\le C (k+1)^{-1-\mu -\alpha /2}\int_0^\infty |f(t)|e^{-t/2}
t^{(\alpha +1)/2+\mu }\, dt, \quad \frac{\gamma +1}{p} \le
\frac{\alpha +1}{2}
-\frac{1}{3p} +\frac{2}{3}  .$$
By Lemma 1, 5th case, in \cite{analy}
$$ \nm{t^{-\mu -1/2} {\cal L}_k^\alpha (t)}_\infty
\le C (k+1)^{\mu +1/2}$$
so that
$$ I\le C (k+1)^{\mu -1/2-\alpha /2}\int_0^\infty |f(t)|e^{-t/2}
t^{(\alpha +1)/2+\mu }\, dt, \quad \quad \quad \quad \quad \quad $$
$$ \le C (k+1)^{-1-\mu -\alpha /2}\int_0^\infty |f(t)|e^{-t/2}
t^{(\alpha +1)/2+\mu }\, dt, \quad \frac{\gamma +1}{p} >
\frac{\alpha +1}{2}
-\frac{1}{3p} +\frac{2}{3}  ,$$
provided that $\mu -(\alpha +1)/2 \le -1-\mu -\alpha /2$ which is
equivalent to
 $\mu \le
-1/4$ or $\gamma \le 3p/4-1/2+\alpha p/2$. But this is no further
restriction
since for $\alpha \le 1/2$ there holds $\alpha +p-1\le
3p/4-1/2+\alpha p/2$.
Summarizing,  for $-1<\alpha \le 1/2,\; \gamma \le \alpha +p-1$ and
$\mu =
(\gamma /p-(\alpha +1)/2)/2(1/p-1/2) $ we have that
\begin{equation}\label{l1est}
\sup_k|(k+1)^{1+\mu +\alpha /2}\Delta_2 {\hat f}_\alpha (k)|
\le C \int_0^\infty |f(t)|e^{-t/2} t^{(\alpha +1)/2+\mu }\, dt.
\end{equation}
 
\medskip \noindent
Now consider  the case $\alpha >1/2$. Then, by formula (2.12) in
\cite{cohen},
(\ref{II}) is obviously true when
$(\gamma +1)/p\le \alpha /2+1 +(1/p-1/2)/2$. Again, the application
of Lemma
1 in \cite{analy} requires $\gamma \le \alpha +p-1$, which for
$\alpha
>1/2$ is less restrictive than $(\gamma +1)/p\le \alpha /2+1
+(1/p-1/2)/2$.
Its 4th case now leads to
$$ I\le C (k+1)^{-11/6-\mu -\alpha /2}\int_0^\infty |f(t)|e^{-t/2}
t^{(\alpha +1)/2+\mu }\, dt, \quad \frac{\gamma +1}{p} \le
\frac{\alpha +1}{2}
-\frac{1}{3p} +\frac{2}{3} , $$
and its 5th case to
$$ I\le C (k+1)^{\mu -1/2-\alpha /2}\int_0^\infty |f(t)|e^{-t/2}
t^{(\alpha +1)/2+\mu }\, dt, \quad (\gamma +1)/p > \frac{\alpha
+1}{2}
-\frac{1}{3p} +\frac{2}{3}  .$$
But $\mu -1/2-\alpha /2 \le -\mu -1-\alpha /2$ if
$(\gamma +1)/p\le \alpha /2+1 +(1/p-1/2)/2$; so that, summarizing,
(\ref{l1est})
also holds under this restriction for $\alpha >1/2$.
 
\medskip \noindent
Now an application of the Stein and Weiss interpolation theorem
(see \cite{sw})
with $Tf=\{ Tf(k)\} $and $ Tf(k)= \sqrt{A_k^{\alpha +1}}
\Delta_2 {\hat f}_\alpha (k)$ gives the assertion of Theorem 1.1 in
the case
$a =0$.
 
\medskip \noindent
If $a \neq 0$ then by (\ref{coeff}), the definition of
$\Delta_2\Delta^a$,
 and by (\ref{gebr})
$$ \Delta _2 \Delta^{a}{\hat f}_\alpha (k)
=C\{\Delta {\hat f}_{\alpha +a} (k)+\Delta{\hat f}_{\alpha +a}(k+1)
\}
=C \Delta _2 {\hat f}_{\alpha +a}(k),$$
since already the condition $\gamma  < \alpha +a+1$ (which implies
no new
restriction) gives absolute
convergence of the infinite sum and integral involved (see the
formula
following (9) in \cite{laguerre}) and Fubini's Theorem can be
applied.
Hence all the previous estimates remain valid when $\alpha$ is
replaced by
$\alpha +a$.
 
\bigskip
\section{ A variant for integrable functions}
 
Theorem 1.1 gives a necessary condition for a sequence $\{ f_k\}$
to generate with respect to $L_k^\alpha$ an $L^1_{w(\gamma
)}$--function. But
this condition is hardly comparable with the following sufficient
one which is
a slight modification of Lemma 2.2 in \cite{laguerre}.
\thm{
Let $\alpha >-1$ and $\delta > 2\gamma -\alpha +1/2\ge 0 $. If
$\{f_k \}$ is a bounded sequence with $\lim_{k\to \infty }f_k =0$
and
$$ \sum_{k=0}^\infty (k+1)^{\delta +\alpha -\gamma}
               |\Delta ^{\delta +1} f_k| \le K_{\{f_k \} } ,$$
then there exists a function $f\in L^1_{w(\gamma )}$ with ${\hat
f}_{\alpha
}(k)=f_k$ for all $k\in {\bf N}_0$ and
$$ \nmgk{f} \le C\, K_{\{f_k \} }$$
for some constant $C$ independent of the sequence $\{f_k \}$.
}
The proof follows along the lines of Lemma 2.2 in \cite{laguerre}
since the
norm of the Ces\`aro kernel
$$ \chi _n^{\alpha ,\delta }(x)=
(A_n^\delta \Gamma (\alpha +1))^{-1}
\sum _{k=0}^n A_{n-k}^\delta L_k^\alpha (x) =
(A_n^\delta \Gamma (\alpha +1))^{-1}L_n^{\alpha +\delta +1}(x) $$
can be estimated with the aid of Lemma 1 in \cite{analy} by
$$ \nm{\chi _k^{\alpha ,\delta }}_{L^1_{w(\gamma )}}
\le C (k+1)^{\alpha - \gamma },\quad  \delta > 2\gamma -\alpha +1/2
$$
 
\medskip \noindent
The variant of Theorem 1.1 in the case $p=1$ is
\thm{
If $\alpha >-1$ and $\gamma >\max \{ -1/3,\, \alpha/2 -1/6\} $,
then
$$\sum_{k=0}^\infty (k+1)^{\gamma -2/3}|\Delta^{2\gamma -\alpha
+1/3}
{\hat f}_\alpha (k)| \le C \nmgk{f}.$$
}
A comparison of the sufficient condition and the necessary one
nicely shows
where the $L^1_{w(\gamma )}$--functions live; in particular we
see that the ``smoothness'' gap (the difference of the orders of
the difference
operators) is just greater than  $7/6$. It is clear that Theorem
3.2
can be modified by using the $\Delta_2$--operator. Theorem 3.2 does
not follow
from the $p=1$ case of Lemma 2.1 in \cite{laguerre} since that
estimate would
lead to the divergent sum $\sum _{k=0}^\infty (k+1)^{-1} \nmgk{f}$.
 
\medskip \noindent
\pf{
By formula (\ref{gebr}) we have
$$ \Delta^{2\gamma -\alpha+1/3}{\hat f}_\alpha (k)
=C \int_0^\infty f(t)R_k^{2\gamma +1/3}(t)
t^{2\gamma +1/3}e^{-t}\, dt \quad \quad \quad \quad $$
$$\quad \quad \quad \quad \quad \quad
=C (k+1)^{-\gamma -1/6}\int_0^\infty f(t){\cal L}_k^{2\gamma
+1/3}(t)
t^{\gamma +1/6}e^{-t/2}\, dt $$
and hence
$$\sum_{k=0}^\infty (k+1)^{\gamma -2/3}|\Delta^{2\gamma
-\alpha+1/3}
{\hat f}_\alpha (k)|
\le C \int_0^\infty |f(t)|\sum_{k=0}^\infty (k+1)^{-5/6}|t^{1/6}
{\cal L}_k^{2\gamma +1/3}(t)|t^{\gamma}e^{-t/2}\, dt$$
if the right hand side converges. To show this
we discuss for $j\in {\bf Z}$
$$ \sup _{2^{j}\le t\le 2^{j+1}} \sum_{k=0}^\infty
(k+1)^{-5/6}|t^{1/6}
{\cal L}_k^{2\gamma +1/3}(t)|$$
and prove that this quantity is uniformly bounded in $j$, whence
the
assertion.
 
\medskip \noindent
First consider those $j\ge 0$ for which there exists a nonnegative
integer $n$
such that $0\le k\le 2^n$ implies $3\nu /2\, :=3(2k+2\gamma +4/3)
\le 2^j$
but such that this inequality fails to hold
for $k\ge 2^{n+1}$; the latter assumption
in particular implies that essentially $\nu /2 \ge 2^{j+1}$ for
$k\ge 2^{n+4}$.
Since
$\nm{t^{1/6}{\cal L}_k^{2\gamma +1/3}(t)}_\infty \le C
(k+1)^{-1/6}$
by Lemma 1 in \cite{analy}, we
 obviously have
\begin{equation}\label{j}
\sum_{k=0}^\infty (k+1)^{-5/6}|t^{1/6}{\cal L}_k^{2\gamma +1/3}(t)|
\le \left( \sum_{k=0}^{2^n} +\sum_{k=2^{n+4}}^\infty \right) \ldots
+ O(1).
\end{equation}
For $k=0,\dots ,2^n$  we can now apply the fourth case of formula
(2.5) in
\cite{analy} to obtain $|t^{1/6}{\cal L}_k^{2\gamma +1/3}(t)|\le C
e^{-\mu 2^j
}$ for some positive constant $\mu$ and the first sum on the right
hand side
of (\ref{j}) is
bounded uniformly in $j$. In consequence of the choice of $n$ the
second
case of formula (2.5) in \cite{analy} can be used for $k\ge
2^{n+4}$,
giving
$$\sum_{k=2^{n+4}}^\infty (k+1)^{-5/6}|t^{1/6}{\cal L}_k^{2\gamma
+1/3}(t)|
 \le C t^{-1/12} \sum_{k=2^{n+4}}^\infty (k+1)^{-13/12} = O(1)$$
since $2^j\le t\le 2^{j+1}$ and $j$ and $n$ are comparable.\\
Now consider the remaining $j$'s: We have to split up the sum
$\sum_{k=0}^\infty \dots $ into two parts, one where $k$ is such
that $2^j\nu
\ge 1$ (this contribution has just been seen to be uniformly
bounded in $j$),
the other where $k$  is such that $2^j\nu \le 1$.
To deal with the last case choose again $n$ to be the greatest
integer
such that $2^{n+2}+4\gamma +8/3 \le 2^{-j}$; this time, $n$ and
$-j$ are
comparable and we obtain by the first case of (2.5) in \cite{analy}
$$\sum_{k=0}^{2^n} (k+1)^{-5/6}|t^{1/6}{\cal L}_k^{2\gamma
+1/3}(t)|
\le C t^{\gamma +1/3}\sum_{k=0}^{2^n}(k+1)^{\gamma -2/3} = O(1)$$
if $2^j\le t\le 2^{j+1},\; \gamma >-1/3$, which completes the
proof.
}

\end{document}